\documentclass[a4paper,twoside,11pt]{amsart}
\usepackage[latin1]{inputenc}
\usepackage[T1]{fontenc}
\usepackage{amssymb,latexsym,array}
\renewcommand{\and}{\unskip{ }et \ignorespaces}

\theoremstyle{plain}
\newtheorem{thm}{\bfseries Théorème}[section]

\newtheorem{prop}[thm]{\bfseries Proposition}

\theoremstyle{remark}

\numberwithin{equation}{section}
\DeclareMathSymbol{\Z}{\mathalpha}{AMSb}{"5A} 
\DeclareMathSymbol{\PP}{\mathalpha}{AMSb}{"50} 
\DeclareMathSymbol{\Q}{\mathalpha}{AMSb}{"51}
\DeclareMathSymbol{\N}{\mathalpha}{AMSb}{"4E}
\DeclareMathSymbol{\R}{\mathalpha}{AMSb}{"52}

\renewcommand{\le}{\leqslant}
\renewcommand{\ge}{\geqslant}

\newcommand{\la}{\lambda}

\newcommand{\ga}{\gamma}

\begin{document}

\title{Quelques calculs de la cohomologie de $GL_N(\Z)$  et  de
 la K-théorie de $\Z$}
\author{Philippe Elbaz-Vincent}
\address{GTA, UMR CNRS 5030, CC51,
Universit\'e Montpellier II, 34095 Montpellier Cedex 5, France.}
\email{pev@math.univ-montp2.fr}

\author{Herbert Gangl}
\address{MPI f\"ur Mathematik Bonn, Vivatsgasse 7, D-53111 Bonn, Deutschland}
\email{herbert@mpim-bonn.mpg.de}
\author{Christophe Soul\'e}
\address{C.N.R.S. et I.H.\'E.S., 35 Route de Chartres, 91440 Bures-sur-Yvette, 
France.}
\email{soule@ihes.fr}
\maketitle
\markboth{{\sc Ph. Elbaz-Vincent, H. Gangl et C. Soulé}}{{\sc Quelques calculs de 
la cohomologie de $GL_N(\Z)$ et de la K-théorie de $\Z$}}
{\small

{{\bf R\'esum\'e} - Pour $N =5$ et $N=6$, nous calculons le complexe
cellulaire d\'efini par Vorono\"i \`a partir des formes quadratiques 
r\'eelles de dimension $N$.
Nous en d\'eduisons l'homologie de  $GL_N(\Z)$ \`a coefficients
triviaux, \`a de petits nombres premiers pr\`es.
Nous montrons aussi que $K_5(\Z) = \Z$ et que $K_6(\Z)$ n'a que de la
 $3$-torsion.}\\[5pt]

\centerline{\normalsize\bf Some computations of the cohomology of
 $GL_N(\Z)$ and the  K-theory of $\Z$}\quad\\[3pt]
\indent {{\bf Abstract} - For $N =5$ and $N=6$, we compute  the Vorono\"i cell 
complex attached to real
$N$-dimensional quadratic forms, 
and we obtain  the homology of $GL_N(\Z)$ with trivial 
coefficients, up to small primes.
We also prove that  $K_5(\Z) = \Z$ and $K_6(\Z)$ has only $3$-torsion.}\\[5pt]


\footnoterule

\section{La th\'eorie de Vorono\"i}
Soit $N\ge 2$ un entier. Notons $C_N$ l'espace des formes quadratiques définies 
positives réelles de rang $N$. \'Etant donnée $h\in C_N$, 
 les vecteurs minimaux de $h$, c'est-\`a-dire les vecteurs $v$ non nuls de 
 $\Z^N$ tels que $h(v)$ soit minimal,
 forment un ensemble fini, not\'e $m(h)$.
 Une forme $h\in C_N$ est dite {\it parfaite} si elle 
 est entièrement caractérisée 
par son minimum sur $\Z^N-\{0\}$ et par l'ensemble $m(h)$.
D\'esignons par $\Gamma$ le groupe  $GL_N(\Z)$ ou $SL_N(\Z)$.
Vorono\"i a démontré 
\cite{voronoi}~(Thm.,\,p.110) que, modulo l'action de $\Gamma$
et la multiplication par les r\'eels positifs, il n'y a 
qu'un nombre fini de formes parfaites.

Notons $C^*_N$ l'espace des formes quadratiques 
positives réelles sur $\R^N$ 
dont le noyau est engendré par un sous-espace vectoriel propre de $\Q^N$. Soient $X^*_N$ le 
quotient de $C^*_N$ par les homothéties positives,  $\pi : C^*_N \to X^*_N$ 
l'application quotient, $X_N = \pi(C_N)$ et $\partial X^*_N = X^*_N - X_N$.
 Le groupe $\Gamma$ agit sur $C^*_N$, et sur $X^*_N$, par la formule
\[h\cdot\ga = \ga^t h \ga,\quad \ga\in \Gamma, \quad h\in C_N^*\, ,\]
où $\ga^t$ est la transposée de $\ga$.\\

A tout vecteur $v\in \Z^N-\{0\}$ on peut  associer une forme $\widehat{v}\in C_N^*$, 
définie par $\widehat{v}(x)=(v|x)^2$. \'Etant donné un sous-ensemble 
fini $B$  de 
$\Z^N-\{0\}$, {\it l'enveloppe convexe} de $B$ est le sous-ensemble de
$X_N^*$
 image par $\pi$ du 
sous-ensemble
\[\Big\{\sum_j \la_j \widehat{v}_j\, ,\, v_j\in B, \, \la_j\ge 0\Big\}\]
de $C_N^*$. Lorsque $h$ est une forme parfaite, nous désignons par 
$\sigma(h)\subset X_N^*$ l'enveloppe convexe de l'ensemble $m(h)$
de ses vecteurs minimaux. 
Voronoï a montré \cite{voronoi}~(\S\S8-15) que les cellules $\sigma(h)$ et leurs 
intersections, quand $h$ parcourt l'ensemble des formes parfaites,
définissent une décomposition cellulaire de $X^*_N$, compatible avec 
l'action de $\Gamma$. Nous munissons  $X^*_N$ de
 la CW-structure correspondante. 
 Si $\tau$ est une cellule (fermée) de $X_N^*$ et si $h$ est une forme parfaite 
telle que $\tau \subset \sigma (h)$, on note $m(\tau)$ l'ensemble des vecteurs 
$v$ de $m(h)$ tels que $\widehat{v}$ soit dans $\tau$. La cellule $\tau$ est 
l'enveloppe convexe de $m(\tau)$ et l'on a $m(\tau) \cap m(\tau') = m(\tau 
\cap \tau')$.
 
\section{Calculs explicites}
\subsection{}
Notons $\Sigma_n$, $0 \le n \le \dim (X_N^*) = d(N)= N(N+1)/2 - 1$, un ensemble de repr\'esentants, 
modulo l'action de $\Gamma$, des 
cellules $\sigma$ de dimension $n$ dans $X_N^*$ qui rencontrent $X_N$ et telles 
qu'aucun élément du stabilisateur de $\sigma$ dans $\Gamma$ 
ne 
change  
l'orientation de $\sigma$.
 Pour $N \le 6$ nous avons déterminé un tel 
ensemble $\Sigma_n$. En particulier on a le résultat suivant~:

\begin{prop} Le cardinal de $\Sigma_n$ est zéro sauf dans les cas indiqués dans le 
tableau ci-dessous~:
\[
\setlength{\extrarowheight}{4pt}
\begin{array}{|c|c|c|c|c|c|c|c|c|c|c|c|c|c|c|c|}
\hline
 n  & 6 & 7 & 8 & 9 & 10 & 11 & 12 & 13& 14 & 15 & 16 & 17 & 18 & 19 & 20 \\
\hline
GL_5(\Z) & 0 & 0 & 1 & 7 & 6 & 1& 0 & 2 & 3 \\
\hline
GL_6(\Z) & 0 &0 & 0 & 3 & 46 & 163 & 340 & 544 & 636 & 469 & 200 & 49 & 5 & 0 & 0\\
\hline
SL_6(\Z) & 3 & 10 & 18 & 43 & 169 & 460 & 815 & 1132 & 1270 
& 970 & 434 & 114 & 27 & 14 & 7\\
\noalign{\hrule height 0.4pt} 
\end{array}
\]
\vskip2pt

\end{prop}

\subsection{}

Pour tout entier $m > 1$ on note ${\mathcal S}_m$ la classe (de Serre) des groupes 
abéliens finis $A$ tels que tout nombre premier $p$ divisant l'ordre de $A$ 
vérifie $p \le m$. Si $\gamma \in \Gamma$ est d'ordre premier $p$ on sait que $p 
\le N+1$. Il en résulte que l'action de $\Gamma$ sur $X_N^*$ permet de définir un 
complexe $V = (V_n , d_n)$ tel que $V_n$ soit isomorphe au module libre engendré 
par $\Sigma_n$, $n \ge 0$, et que l'homologie de $V$ co{\"\i}ncide, modulo 
${\mathcal S}_{N+1}$, avec l'homologie $\Gamma$-équivariante de la paire $(X_N^* , 
\partial X_N^*)$ (\cite{B} VII Prop.~8.1, 
\cite{soule-3torsionK4} Prop.~2.2)~:
\begin{equation*}
H_n (V) = H_n^{\Gamma} (X_N^* , \partial X_N^* ; \Z ) \qquad ({\rm mod} \, {\mathcal 
S}_{N+1}) \, .
\end{equation*}
\vskip 10pt
\begin{prop}
\begin{itemize}
\item[i)] Si $\Gamma = {\rm GL}_5 (\Z)$, on a, modulo ${\mathcal S}_5$,
\[H_n(V,\Z)=\begin{cases} \Z &{\rm si }\quad n=9,14,\\
0 & {\rm sinon}.
\end{cases}\]

\item[ii)] Si $\Gamma = {\rm GL}_6 (\Z)$, on a, modulo ${\mathcal S}_7$,

\[H_n(V,\Z)=\begin{cases} 
 \Z &{\rm si }\quad n=10, 11, 15,\\
0 & {\rm sinon}.
\end{cases}\]
\item[iii)] Si $\Gamma = {\rm SL}_6 (\Z)$, on a, modulo ${\mathcal S}_7$,

\[H_n(V,\Z)=\begin{cases} \Z^2 &{\rm si }\quad n=15,\\
 \Z &{\rm si }\quad n=10, 11, 12, 20,\\
0 & {\rm sinon}.
\end{cases}\]
\end{itemize}
\end{prop}

\subsection{}

La preuve des Propositions 2.1 et 2.2 utilise la classification des formes parfaites. 
Si $N \le 7$, le travail de Jaquet \cite{jaquet} fournit un ensemble ${\mathcal 
P}$ de représentants des formes parfaites  de rang $N$ et, si $h \in {\mathcal 
P}$, la liste $m(h)$ des vecteurs minimaux de $h$ et une liste des formes 
``voisines'' de $h$ modulo l'action de son stabilisateur $\Gamma_h$, ainsi que 
leurs représentants dans ${\mathcal P}$. On en d\'eduit une liste 
des faces de $\sigma (h)$ modulo $\Gamma_h$.

Pour obtenir un ensemble $\Sigma_n$ comme dans 2.1,
 on procède comme suit, à l'aide 
d'un ordinateur. Pour chaque $h \in {\mathcal P}$ on calcule les 
éléments de $\Gamma_h$, 
 et donc la liste ${\mathcal F}_{1,h}$ de toutes ses faces $\tau$ de 
codimension 1 (i.e. les ensembles $m(\tau)$ correspondants).
Par récurrence sur l'entier $n$,
 $1 \le n \le d(N)$, on définit 
un ensemble ${\mathcal F}_{n,h}$ de cellules de codimension $n$ dans $X_N^*$, puis un 
système de représentants ${\mathcal C}_{n,h} \subset {\mathcal F}_{n,h}$ pour 
l'action de $\Gamma$ sur ${\mathcal F}_{n,h}$ (cf. 4.3). Par définition, ${\mathcal 
F}_{n+1,h}$ est l'ensemble des cellules $\varphi \cap \tau$, $\varphi \in {\mathcal 
F}_{n,h}$, $\tau \in {\mathcal C}_{n,h}$, qui sont de codimension $n+1$ dans 
$\sigma (h)$.
L'ensemble $\Sigma_n$ est alors obtenu en choisissant des représentants modulo 
$\Gamma$ de la réunion des ${\mathcal C}_{d(N) - n,h}$, $h \in {\mathcal P}$.

\subsection{}

Pour calculer la diff\'erentielle de $V$ on procède comme suit. Pour chaque 
cellule $\sigma$ de $\Sigma_n$, $n \ge 0$, on choisit un ordre sur l'ensemble 
$m(\sigma)$ de ses vecteurs minimaux. Pour toute cellule $\tau' \subset \sigma$ 
cet ordre fournit un ordre sur $m(\tau')$ et donc une orientation du sous-espace 
vectoriel réel $\R \, \langle \tau' \rangle$ engendré par $m(\tau')$ 
dans l'espace vectoriel
des matrices symétriques réelles. Si $\tau'$ est une face de $\sigma$ on obtient 
une base orientée $B$ de $\R \, \langle \sigma \rangle$ en adjoignant à la suite 
d'une base positive de $\R \, \langle \tau' \rangle$ l'élément $\widehat v$, où $v$ 
est l'élément minimal 
de $m(\sigma) - m(\tau')$
(pour l'ordre choisi ci-dessus). 
Notons $\epsilon (\tau' , \sigma) 
= \pm 1$ l'orientation de $B$ dans $\R \, \langle \sigma \rangle$.

Si $\tau \in \Sigma_{n-1}$ est équivalente à la face $\tau' = \tau \cdot \gamma$ 
de $\sigma \in \Sigma_n$, on pose $\eta (\tau , \tau') = 1$ (resp. $\eta 
(\tau,\tau') = -1$) selon que $\gamma$ est compatible (ou non) aux orientations 
choisies de $\R \, \langle \tau \rangle$ et $\R \, \langle \tau' \rangle$. On a alors
$$
d_n (\sigma) = \sum_{\tau \in \Sigma_{n-1}} \sum_{\tau'} \eta (\tau,\tau') \, 
\epsilon (\tau' , \sigma) \, \tau \, ,
$$
où $\tau'$ parcourt l'ensemble des faces de $\sigma$ équivalentes à $\tau$. On 
notera que l'identité $d_{n+1} \circ d_n = 0$, $n \ge 0$, est un bon test pour 
vérifier les calculs effectués par l'ordinateur.

\section{Cohomologie des groupes modulaires}

Si ${\rm St}_N$ est le module de Steinberg de ${\rm SL}_{N,\Q}$ on a
\begin{equation}
H_n^{\Gamma} (X_N^* , \partial X_N^* ; \Z) = H_{n-N+1} (\Gamma , {\rm St}_N)
\end{equation}
(cf. \cite{soule-3torsionK4}). Le théorème de dualité de Borel et Serre
 \cite{BS} dit que
$$
H_m (\Gamma , {\rm St}_N) = H^{d-m} (\Gamma , \widetilde{\Z}) \qquad ({\rm mod} \, {\mathcal 
S}_{N+1}) \, ,
$$
où $d= N(N-1)/2$ est la dimension cohomologique virtuelle
 de $\Gamma$ et $\widetilde{\Z}$ 
est le $\Gamma$-module trivial $\Z$ sauf quand $\Gamma = {\rm GL}_N (\Z)$ et $N$ est 
pair, auquel cas $\gamma \in \Gamma$ agit par multiplication par $\det (\gamma)$ 
sur $\widetilde{\Z} = \Z$. Enfin, le lemme de Shapiro montre que, modulo 
${\mathcal S}_2$, 
$$
H^m ({\rm SL}_N (\Z) , \Z) = H^m ({\rm GL}_N (\Z) , \Z) \oplus H^m ({\rm GL}_N 
(\Z) , \widetilde{\Z}) \, .
$$
La Proposition 2.2 implique donc le résultat suivant~:

\begin{thm}
\begin{itemize}
\item[i)] Modulo ${\mathcal S}_5$ on a
$$
H^m ({\rm GL}_5 (\Z) , \Z) = \begin{cases} \Z & {\rm si} \quad m = 0, 5,\ \\
0 & {\rm sinon}.
\end{cases}
$$
\item[ii)] Modulo ${\mathcal S}_7$ on a
$$
H^m ({\rm GL}_6 (\Z) , \Z) = \begin{cases} \Z &{\rm si }\quad m= 0, 5, 8,\\
0 & {\rm sinon},
\end{cases}
$$
et
$$
H^m ({\rm SL}_6 (\Z) , \Z) = \begin{cases} \Z^2 &{\rm si }\quad m=5,\\
\Z &{\rm si }\quad m=0, 8, 9, 10,\ \\
0 & {\rm sinon}.
\end{cases}
$$
\end{itemize}
\end{thm}

\section{K-th\'eorie des entiers}

\subsection{}

Rappelons que $K_1 (\Z) = K_2 (\Z) = \Z / 2$, $K_3 (\Z) = \Z / 48$ et $K_4 (\Z) = 
0$ \cite{R}.
\vskip 10pt
\begin{thm}
On a $K_5 (\Z) = \Z$. Par ailleurs, l'ordre de $K_6 (\Z)$ est une puissance de 3.
\end{thm}
\vskip 10pt

Notons $Q$ (resp. $Q_N$) la catégorie définie par Quillen \cite{Q1} à partir des 
$\Z$-modules libres de rang fini (resp. de rang au plus $N$). Si $BQ$ est le 
classifiant de $Q$, on a
$$
K_m (\Z) = \pi_{m+1} \, BQ \, ,
$$
et l'on dispose de suites exactes \cite{Q2}
\begin{equation}
\cdots \rightarrow H_m (BQ_{N-1} , \Z) \rightarrow H_m (B Q_N , \Z) \rightarrow 
H_{m-N} ({\rm GL}_N (\Z) , {\rm St}_N)
 \rightarrow H_{m-1} (B Q_{N-1} , \Z) \rightarrow \cdots 
\end{equation}

On sait aussi que $H_0 ({\rm GL}_N (\Z) , {\rm St}_N) = 0$ si $N \ge 1$.

\begin{prop}
\begin{itemize}
\item[i)] Modulo ${\mathcal S}_2$ on a
$$
H_3 ({\rm GL}_3 (\Z) , {\rm St}_3) = \Z
$$
et
$$
 H_1 ({\rm GL}_5 (\Z) , {\rm St}_5) = H_2 ({\rm GL}_4 (\Z) , {\rm St}_4) =
 H_4 ({\rm GL}_2 (\Z) , {\rm St}_2) = 0 \, .
$$
\item[ii)] Modulo ${\mathcal S}_3$ on a
$$
H_3 ({\rm GL}_4 (\Z) , {\rm St}_4) = \Z
$$
et
$$
H_1 ({\rm GL}_6 (\Z) , {\rm St}_6) = H_2 ({\rm GL}_5 (\Z) , {\rm St}_5) = H_4 
({\rm GL}_3 (\Z) , {\rm St}_3) = H_5 ({\rm GL}_2 (\Z) , {\rm St}_2) = 0 \, .
$$
\end{itemize}
\end{prop}

Pour d\'emontrer cette proposition on utilise les r\'esultats du calcul amenant \`a
la Proposition 2.1 et les résultats de 
\cite{SL3}, \cite{lee/szczarba},  \cite{soule-torsionKZ}
et \cite{soule-3torsionK4}. 
D'après (3.1) et \cite{B} ou \cite{soule-3torsionK4}, on peut 
calculer les groupes $H_m ({\rm GL}_N (\Z) , {\rm St}_N)$ à l'aide d'une suite 
spectrale dont le terme $E_1$ est une somme de groupes d'homologie des 
stabilisateurs des cellules des $\Sigma_n$. L'analyse 
de ces groupes conduit à la 
Proposition 4.2.

\subsection{}

A l'aide de (4.1) et de \cite{lee/szczarba} on déduit de la Proposition 4.2 que $H_6 (BQ,\Z) = 
\Z$ modulo ${\mathcal S}_2$ et que $H_7 (BQ,\Z) = \Z$ modulo ${\mathcal S}_3$. Par 
ailleurs on d\'emontre, en utilisant \cite{A},  que le noyau du morphisme d'Hurewicz
$$
K_5 (\Z) \rightarrow H_6 (BQ,\Z)
$$
(resp. $K_6 (\Z) \rightarrow H_7 (BQ,\Z)$) est dans ${\mathcal S}_2$ (resp. 
${\mathcal S}_3$). Enfin, on sait que $K_6 (\Z)$ est fini, que $K_5 (\Z)$ 
est de rang un, et qu'aucun d'eux n'a de la 2-torsion \cite{RW}. Le Théorème 4.1 en 
résulte.

\subsection{}

La plupart des programmes ont \'et\'e d\'evelopp\'es en utilisant le logiciel PARI-GP
\cite{parigp}. Nous utilisons aussi les programmes de 
 \cite{souvignier} pour d\'ecider si deux formes de $C_N^*$ sont
 \'equivalentes sous l'action de $\Gamma$. Enfin, le logiciel
GAP \cite{GAP} permet de produire tous les \'el\'ements des groupes finis 
$\Gamma_h$, $h \in {\mathcal P}$, \`a partir de leurs g\'en\'erateurs et 
de calculer certains de leurs groupes d'homologie.

Les calculs ont \'et\'e r\'ealis\'es sur les ordinateurs de l'UMS M\'edicis, 
ceux du MPI Bonn,
et ceux de l'ACI "Arithm\'etique des fonctions L".
}

\bibliographystyle{plain}

\end{document}